\documentclass[3p, 11pt]{elsarticle} 
\usepackage[english]{babel}
\usepackage{natbib}
\usepackage[ruled,linesnumbered]{algorithm2e}
\usepackage{amssymb,amsmath,graphicx}
\usepackage{amsthm}
\usepackage{float}
\usepackage{afterpage}
\usepackage{changepage}
\usepackage{dcolumn}
\usepackage{multirow}
\usepackage{color}
\usepackage[pdftex]{pict2e}
\usepackage{enumitem}
\usepackage{longtable}
\usepackage{indentfirst}
\usepackage{mathtools}
\usepackage{subcaption}
\usepackage{tabularx}
\usepackage{changes}
\usepackage{multicol}
\usepackage{tikz}
\usetikzlibrary{shapes.geometric, arrows}
\usetikzlibrary{matrix}
\usepackage{adjustbox} 

\tikzstyle{startstop} = [rectangle, rounded corners, minimum width=3cm, minimum height=1cm,text centered, draw=black, fill=gray!20]
\tikzstyle{process} = [rectangle, minimum width=3.5cm, minimum height=1cm, text centered, draw=black, fill=blue!10]
\tikzstyle{decision} = [diamond, minimum width=3cm, minimum height=1cm, text centered, draw=black, fill=red!10, inner sep=0pt]
\tikzstyle{arrow} = [thick,->,>=stealth]

\usepackage[figuresright]{rotating}

\usepackage{hyperref}       
\usepackage{url}            
\usepackage{booktabs}       
\usepackage{nicefrac}       
\usepackage{microtype}      
\usepackage{lipsum}
\usepackage{graphicx}       
\usepackage{listings}       
\usepackage{xcolor}         
\usepackage{multirow}       
\usepackage[pdftex]{pict2e}
\usepackage{float}

\SetKwComment{Comment}{\color{green!50!black}// }{}

\SetKwProg{Function}{Function}{}{end}

\allowdisplaybreaks

\definecolor{codegreen}{rgb}{0,0.6,0}
\definecolor{codegray}{rgb}{0.5,0.5,0.5}
\definecolor{codepurple}{rgb}{0.58,0,0.82}
\definecolor{backcolour}{rgb}{0.95,0.95,0.92}

\begin{document}

\begin{frontmatter}




\title{A Hybrid Matheuristic Framework for the Chinese Postman Problem with Load-Dependent Costs}


\author[1]{Thieu Khang Nguyen}
\author[2]{Thu Huong Dang \texorpdfstring{\corref{cor1}}}
\author[1]{Truong-Son Hy \texorpdfstring{\corref{cor1}}}

\cortext[cor1]{Correspondence to thy@uab.edu and t.h.dang@lancaster.ac.uk}
\affiliation[1]{organization={University of Alabama at Birmingham},
postcode={AL 35294},
city={Birmingham},
country={The United States}}
\affiliation[2]{organization={Lancaster University},
postcode={LA1 4YX},
city={Lancaster},
country={The United Kingdom}}

\address{}

\begin{abstract}

The Chinese Postman Problem with load-dependent costs (CPP-LC) arises in real-world logistics and transportation systems where travel costs depend on vehicle load and energy consumption. In this work, we propose a hybrid optimization framework that integrates metaheuristic search with mathematical programming to efficiently solve CPP-LC. The proposed method combines local search procedures with reduced mixed-integer linear programming (MILP) models to balance exploration and intensification. In addition, we develop an Ant Colony Optimization (ACO) algorithm to enhance scalability on large instances. Extensive experiments on benchmark datasets demonstrate that the proposed framework consistently achieves high-quality solutions and outperforms existing approaches in solution quality, while maintaining competitive computational efficiency. These results highlight the effectiveness of hybrid optimization strategies for complex, load-dependent routing problems in practical applications.
Our implementation is publicly available at \url{https://github.com/HySonLab/MatCPP}.

\end{abstract}


\end{frontmatter}


\section{Introduction}
\label{sec:intro}

Efficient transportation and logistics planning are crucial to modern society, driving extensive research on formulations, theories, applications, and algorithms of \textit{Vehicle Routing Problems} (VRPs) and their numerous variants (see the books \cite{GRW08,TV14}). Among these, \textit{Arc Routing Problems} (ARPs) form a special class of VRPs where demands are associated with the edges or arcs of a network rather than its nodes. ARPs arise in a wide range of real-world applications, including postal delivery, meter reading, refuse collection, winter road maintenance, and street sweeping (see the books \cite{CL15a,Dr00} and the surveys \cite{Co21, MP17}). Such problems are increasingly relevant in modern intelligent transportation systems, where energy-aware routing and sustainability considerations play a critical role.

In routing literature, the cost of traveling an edge is typically assumed to depend only on its length. However, in practice, the vehicle's weight at the time of traversal also affects fuel consumption and emissions, thereby influencing travel costs. In fact, fuel consumption and emission models, such as those in \cite{BYS05, DBL11, LTL12}, show that both heavier vehicles and vehicles carrying heavier payloads require more energy to move, leading to higher fuel consumption and, consequently, higher pollutant emissions. Therefore, integrating the \textit{gross vehicle weight} (curb weight plus payloads) into travel cost calculations is crucial, especially given the increasing research interest in addressing environmental concerns in transportation.

In the VRP literature, problems in which the travel cost depends on the vehicle's gross weight are known as load-dependent VRPs. There are several works dealing with these problems (e.g., \cite{xiao2012development,zachariadis2015load,fukasawa2016branch,liu2019constraint,rastani2023large}), but their arc routing counterparts have been explored to a much lesser extent. The only study to date that explores this area of research is by Corber\'an et al. \cite{corberan2018}, who investigated the load-dependent Chinese Postman Problem (CPP-LC). The CPP is a basic ARP, where all edges must be traversed at a minimum cost. For a formal definition and algorithms of the CPP, the reader is referred to the book \cite{Dr00}. Unlike the CPP, which is polynomially solvable, the CPP-LC is generally strongly NP-hard.

Corber\'an et al. \cite{corberan2018} present two mathematical programming formulations for the CPP-LC: one that uses characteristics of the ARP, called the arc routing formulation, and another that transforms the original problem into an equivalent VRP, referred to as the node routing formulation. The node routing formulation is linear, whereas the arc routing formulation is originally non-linear but can be linearized. The experimental results show that, although the node routing formulation involves fewer constraints, the arc routing formulation outperforms it when solved using CPLEX’s MILP solver. Moreover, the findings highlight that even for small instances with fewer than 20 nodes and 48 edges, CPLEX was unable to find an optimal solution in an hour of computation time. To address this, the authors proposed two widely known metaheuristics, \textit{Iterated Local Search} (ILS) and \textit{Variable Neighborhood Search} (VNS). When evaluating local search operators within these metaheuristics, a linear-time dynamic programming approach is proposed to efficiently determine the traversal direction of each edge in a given sequence. The experimental results show that for small-scale instances where the CPLEX solver can find a feasible solution within one hour, both metaheuristics are able to generate better solutions after one hour of running time.

In this paper, we study the CPP-LC. The main contribution of this paper is the design of a heuristic solution for the CPP-LC, which combines the advantages of a heuristic procedure for quickly finding local optima with the advantages of using MILP models to intensify the search in specific regions of the solution space. The combination of heuristic schemes with MILP models was first surveyed by Ball \cite{ball2011heuristics} for combinatorial optimization problems in general, and more recently by Archetti and Speranza \cite{archetti2014survey} for routing problems. We adopt the term ``matheuristic'', as introduced in \cite{boschetti2010combining, boschetti2023matheuristics}. This work can be viewed as a hybrid optimization framework bridging operations research and intelligent computational methods.

Throughout our matheuristic, we adopt the arc routing formulation introduced in \cite{corberan2018} to construct MILP models to intensify the search. To improve computational efficiency, we generate small-size models that can be solved to optimality more quickly. Additionally, while Corber\'an et al. (2018) proposed two well-known metaheuristics for the CPP-LC, we contribute another widely known metaheuristic, \textit{Ant Colony Optimization} (ACO), for the sake of completeness.

We conducted experiments on a set of benchmark instances generated following the approach in \cite{corberan2018}, as well as on artificially generated instances. The results confirm the effectiveness of our matheuristic. However, the performance of ACO varies with the size of the instance. We also conduct a sensitivity analysis to evaluate the impact of the main components of our matheuristic.

The paper is structured as follows: Section \ref{sec:review} describes the CPP-LC, while Section \ref{sec:methodology} presents the matheuristic scheme and its components. The ACO metaheuristic is detailed in Section \ref{sec:aco}. The computational results and a discussion are provided in Section \ref{sec:experiments}. Section \ref{sec:conclusion} concludes the study.
\section{The Chinese Postman Problem with load-dependent costs}
\label{sec:review}

\subsection{Problem description}

The CPP-LC is a variant of the classical CPP in which the cost of traversing an edge depends on both its length and the current load of the vehicle. Let $G = (V, E)$ be an undirected connected graph, where $V = \{1, \dots, n\}$ is the set of vertices and $E = \{e_1, \dots, e_m\}$ is the set of edges. Vertex 1 denotes the depot. Each edge $e \in E$ is associated with a length $d_e \geq 0$ and a demand $q_e \geq 0$.

A vehicle with curb weight $W$ is initially loaded with total demand $Q = \sum_{e \in E} q_e$. Starting from the depot, the vehicle must traverse each edge at least once to service the demands and return to the depot. The first time an edge $e = (i, j)$ is traversed, it is serviced. Let $f_e$ be the vehicle load at node $i$, immediately before traversing edge $e$. Then:
\begin{itemize}
    \item The cost of deadheading $e$ (that is, traversing without service) is $d_e(W + f_e)$.
    \item The cost of service $e$ is $d_e\left(W + f_e - \frac{q_e}{2}\right)$.
\end{itemize}

The objective of the CPP-LC is to minimize the total cost. 

\subsection{MILP formulation}
\label{sec:milp}

We now present the MILP version of the arc routing model introduced by \cite{corberan2018}. Since the total cost of a tour depends on the order in which the required edges are traversed, the route is partitioned into $m$ segments (referred to as ``periods''). Each segment begins with the service of exactly one required edge and may include any number of deadheaded edges before reaching the next required edge. As a result, each segment contains exactly one serviced edge, and any feasible route is composed of $m$ such segments.

The MILP formulation uses the following variables for each edge $(i, j) \in E$ and each segment $1 \leq k \leq m$:

\begin{itemize}
    \item The binary variable $y_{i,j}^{k}$ is equal to $1$ if the edge $(i, j) \in E$ is served at the beginning of the segment $k$ and $0$ otherwise.
    \item The binary variable $x_{i,j}^{k}$ is equal to $1$ if the edge $(i, j) \in E$ is deadheaded during the segment $k$ and $0$ otherwise.
    \item The continuous variable $f_k$ denotes the vehicle load at the beginning of the segment $k$.
    \item The continuous variable $l_e^k$ represents the vehicle load at the beginning of the segment $k$ if the edge $e \in E$ is served in that segment (that is, $f_{k}(y_{i,j}^k + y_{j,i}^k) = l_e^k$).
    \item The continuous variable $r_e^k$ denotes the vehicle load after servicing the required edge in segment $k$, if edge $e = (i,j) \in E$ is deadheaded during the same segment (that is, $f_{k+1}(x_{i,j}^k + x_{j,i}^k) = r_e^k$).
\end{itemize}

The MILP is as follows:

\newcommand{\forallshift}{-6cm} 
\begin{align}
  \min \quad & \sum_{k=1}^{m} \sum_{e=(i, j) \in E}
        d_e \left( (W - q_e/2)(y_{ij}^k + y_{ji}^k)
        + W(x_{ij}^k + x_{ji}^k) \right)
        + \sum_{k=1}^{m} \sum_{e \in E} d_e l_e^k
        + \sum_{k=1}^{m-1} \sum_{e \in E} d_e r_e^k. \label{eq:objective_final} \\
  \text{s.t.} \quad 
  & \sum_{k=1}^{m} (y_{ij}^k + y_{ji}^k) = 1 &&\hspace{\forallshift} \quad \forall e \in E  \label{eq:constraint1} \\
  & \sum_{(i, j) \in E} (y_{ij}^k + y_{ji}^k) = 1 &&\hspace{\forallshift} \quad \forall k \in \{1, \dots, m\} \label{eq:constraint2} \\
  & f_{k+1} = f_k - \sum_{e=(i,j) \in E} q_e (y_{ij}^k + y_{ji}^k) &&\hspace{\forallshift} \quad \forall k \in \{1, \dots, m\} \label{eq:constraint3} \\
  & f_1 = Q, \quad f_{m+1} = 0 \label{eq:constraint4} \\
  & y^{1}(\delta^-(i)) + x^{1}(\delta^-(i)) = y^{2}(\delta^+(i)) + x^{1}(\delta^+(i)) &&\hspace{\forallshift} \quad \forall i \in V \setminus \{1\}  \label{eq:constraint5} \\
  & y^{k}(\delta^-(i)) + x^{k}(\delta^-(i)) = y^{k+1}(\delta^+(i)) + x^{k}(\delta^+(i)) &&\hspace{\forallshift} \quad \forall i \in V, k \in \{2, \dots, m-1\}  \label{eq:constraint6} \\
  & y^{m}(\delta^-(i)) + x^{m}(\delta^-(i)) = x^{1}(\delta^+(i)) &&\hspace{\forallshift} \quad \forall i \in V \setminus \{1\}  \label{eq:constraint7} \\
  & y^{m}(\delta^-(1)) + x^{m}(\delta^-(1)) = y^{1}(\delta^+(1)) = 1 \label{eq:constraint8} \\
  & x_{ij}^k \leq \sum_{l=1}^{m-1} (y_{ij}^l + y_{ji}^l) + y_{ji}^k &&\hspace{\forallshift} \quad \forall (i,j) \in E,  k \in \{1, \dots, m\} \label{eq:constraint9} \\
  & L_k (y_{ij}^k + y_{ji}^k) \le l_e^k \leq f_k + L_k (y_{ij}^k + y_{ji}^k - 1) &&\hspace{\forallshift} \quad \forall (i,j) \in E,  k \in \{1, \dots, m\}\label{eq:constarintExtra1}\\
  & f_k + U_k (y_{ij}^k + y_{ji}^k - 1)  \leq l_e^k \leq U_k (y_{ij}^k + y_{ji}^k) &&\hspace{\forallshift} \quad \forall (i,j) \in E,  k \in \{1, \dots, m\} \\
  & L_{k+1} (x_{ij}^k + x_{ji}^k) \le r_e^k \le f_{k+1} + L_{k+1} (x_{ij}^k + x_{ji}^k - 1) &&\hspace{\forallshift} \quad \forall (i,j) \in E,  k \in \{1, \dots, m-1\} \\
  & f_{k+1} + U_{k+1} (x_{ij}^k + x_{ji}^k - 1) \le r_e^k \leq U_{k+1} (x_{ij}^k + x_{ji}^k) &&\hspace{\forallshift} \quad \forall (i,j) \in E, k \in \{1, \dots, m-1\} \\
 & L_k \leq f_k \leq U_k &&\hspace{\forallshift} \quad \forall k \in \{2, \dots, m\} \label{eq:constraint11} \\
  & x_{ij}^k, x_{ji}^k, y_{ij}^k, y_{ji}^k \in \{0, 1\} &&\hspace{\forallshift} \quad \forall (i, j) \in E,  k \in \{1, \dots, m\} \label{eq:constraint10} 
\end{align}

Constraints \eqref{eq:constraint1}–\eqref{eq:constraint2} ensure that each required edge is served exactly once across all segments, and each segment includes exactly one required edge. The constraint \eqref{eq:constraint3} updates the remaining flow capacity by subtracting the demand served in each segment. The constraint \eqref{eq:constraint4} specifies that the first segment starts with the total demand and the final segment ends with zero remaining flow. Constraints \eqref{eq:constraint5}–\eqref{eq:constraint8} enforce flow conservation at each node across and within segments, where $x^{k}(\delta^-(v))=\sum_{(i, j) \in E}{x_{ji}^k}$, $x^{k}(\delta^+(v))=\sum_{(i, j) \in E}{x_{ij}^k}$, and the same applies for $y^{k}(\delta^-(v))$ and $y^{k}(\delta^+(v))$. The constraint \eqref{eq:constraint9} ensures that an edge can only be deadheaded if it has been previously served. Constraints~\eqref{eq:constarintExtra1}--\eqref{eq:constraint11} define and tighten $l_e^k$, $r_e^k$, and $f_k$, where $L_k = Q - \displaystyle \max_{\substack{S \subseteq E \\ |S| = k}} \sum_{e \in S} q_e$ and $U_k = Q -  \displaystyle \min_{\substack{S \subseteq E \\ |S| = k}} \sum_{e \in S} q_e$. The remaining constraints specify the ranges for all variables.

\subsection{Existing algorithms for CPP-LC}
\label{subse:existingAlgs}
To calculate the total cost of a given sequence, Corber\'an et al. \cite{corberan2018} introduced a dynamic programming to find the best direction to traverse each edge (either forward or backward) that minimizes the overall cost. This method runs in $O(m)$ time. However, this method requires one to know the shortest distances between every pair of vertices in the graph beforehand. The calculation of these shortest paths is done using the Floyd–Warshall algorithm, which takes $O(n^3)$ time.

The dynamic programming formulation is as follows. Let $S = ( e_1, \ldots, e_m)$ be a sequence of edges, and let $S_k$ denote the subsequence of $S$ starting from the edge $e_k = \{i_k, j_k\}$. For each $k \in \{1, \ldots, m\}$ and direction $d \in \{1, 2\}$, let $z_d(S_k)$ denote the minimum cost of service $S_k$, assuming that the edge $e_{k}$ is traversed in direction $d$. Here, $d = 1$ indicates traversing $e_k$ from $i_k$ to $j_k$, and $d = 2$ from $j_k$ to $i_k$. We define a dummy edge $e_{m+1} = \{0, 0\}$ with zero cost and load.

For $1 \leq k \leq m$, the values $z_{k,1}$ and $z_{k,2}$ are computed recursively as follows:

\begin{equation*}
\begin{aligned}
z_1(S_k) &= d_{e_k} \left(W + f_k - \frac{q_{e_k}}{2}\right) + 
\min \Big( D_{j_{k}, i_{k+1}} (W + f_{k+1}) + z_1(S_{k+1}),\,
           D_{j_{k}, j_{k+1}} (W + f_{k+1}) + z_2(S_{k+1}) \Big), \\
z_2(S_k) &= d_{e_k} \left(W + f_k - \frac{q_{e_k}}{2}\right) + 
\min \Big( D_{i_{k}, i_{k+1}} (W + f_{k+1}) + z_1(S_{k+1}),\,
           D_{i_{k}, j_{k+1}} (W + f_{k+1}) + z_2(S_{k+1}) \Big).
\end{aligned}
\label{eq:fi_group}
\end{equation*}

Here, $D_{ij}$ is the shortest distance from $i\in V$ to $j\in V$, $z_1(S_{m+1}) = WD_{j_m,0}$, $z_2(S_{m+1}) = WD_{i_m,0}$, and $f_k$ denotes the total demand of $S_k$. The minimum cost of servicing the required edges in a given sequence $S$ is $z(S) = \min \Big(D_{1, i_1}(W+Q) + z_1(S_1), D_{1, j_1}(W+Q) + z_2(S_1) \Big)$. In short, $z(S)$ is referred to as the \textit{total cost} of the sequence $S$.

The authors then proposed the \textit{Greedy Constructive Heuristic} (GCH) for the CPP-LC. The GCH begins by sorting all edges in non-increasing order based on the product of their demand and distance. It then constructs the solution greedily as follows. Starting with an empty sequence, each edge from the sorted list is considered in turn and inserted into the position within the current partial sequence that yields the smallest increase in total cost. To evaluate each potential insertion position, a dynamic programming approach is employed. This procedure is repeated until all edges have been inserted.

For large-scale CPP-LC instances, the authors proposed the ILS and VNS metaheuristics. Both methods start from a solution generated by the GCH and perform up to $k_{\max}$ iterations. In each iteration, the current best solution is perturbed by applying $\frac{m}{5}$ random edge exchanges, followed by a local search. ILS simultaneously evaluates the 2-Exchange, 1-Opt, and 2-Opt operators and applies the best. In contrast, VNS applies these operators sequentially in that specified order and accepts the first move that yields an improvement. If the resulting solution improves on the best solution, it is accepted. For full details of these algorithms, we refer the reader to \cite{corberan2018}.

\section{Methodology}
\label{sec:methodology}

In this section, we describe a matheuristic for the CPP-LC that we call MaLD. Before describing the MaLD scheme and its components, we first introduce some notation.

\subsection{Notation}

A CPP-LC route consists of a sequence of required edges traversed in order, the traversal direction for each required edge, and the deadheading (non-servicing) paths between consecutive edges, as well as to and from the depot. However, storing these full details is often unnecessary in MaLD, as it increases memory requirements and computational complexity.

We use the term \textit{sequence} to refer to an ordered list of edges that must be traversed, without specifying details such as traversal directions or deadheading paths. In contrast, the term \textit{route} refers to the complete route, starting and ending at the depot, where these details are known. For example, a sequence $S = (e_1, e_2, \dots, e_m)$ lists the edges $e_1, e_2, \dots, e_m$ to be visited sequentially, while its corresponding route is the least-cost path that starts from the depot, deadheads to the starting point of $e_1$, services $e_1$, then $e_2$, and so on, until $e_m$, before deadheading back to the depot. To compute the route for a given sequence, we use dynamic programming to determine the traversal direction for each edge. The shortest paths between consecutive edges, as well as between the depot and the last edges, are precomputed before executing MaLD.

\subsection{MaLD scheme}
MaLD is a matheuristic that combines a local search procedure with the exact solution of MILP models.

The algorithm begins by generating an initial solution using a greedy heuristic detailed in Subsection \ref{subse:existingAlgs}. The initial solution is iteratively refined through two main procedures, local search and intensification, until a predefined stopping criterion is met. The local search procedure quickly explores the neighborhood of the current solution to find improved solutions. The intensification phase uses the MILP models to intensify the search to further improve the solution. After the local search procedure, an edge reordering step is performed to prepare the solution for input into the MILP model. A general scheme of MaLD is provided in Algorithm~\ref{alg:matheuristic}.

We now describe in detail the different components that compose MaLD.

\begin{algorithm}[tbh]
\caption{MaLD scheme}
\label{alg:matheuristic}
\KwIn{
    Graph $G = (V, E)$, number of outer iterations $k_{\max}^1$, 
    maximum number of intensification attempts per outer iteration $k_{\max}^2$}
\KwOut{Sequence $S$}

Generate an initial sequence $S$ using a greedy constructive heuristic\;
\For{$i = 1$, \ldots, $k_{\max}^1$}{
    Apply the \texttt{local search procedure} to $S$ (Algorithm~\ref{alg:m-shift})\;
    Apply the \texttt{reordering procedure} to $S$ (Algorithm~\ref{alg:preprocess})\;
    
    \For{$j = 1$, \ldots, $k_{\max}^2$}{
        Apply the \texttt{intensification procedure} to $S$\;
        \If{an improvement is found}{
            Update sequence $S$\;
            \textbf{break}\;
        }
    }
}
\end{algorithm}

\subsection{Local search procedure}

Once the initial solution has been obtained, a local search procedure is carried out. Algorithm~\ref{alg:m-shift} gives details about this procedure.

This procedure potentially improves the current solution by applying a general shift move. Specifically, a set $T$ of $\theta$ edges is randomly selected from the current sequence and removed. Then each edge in $T$ is reinserted, one at a time, in the position that minimizes total cost. The cost of each possible reinsertion is efficiently evaluated using a dynamic programming algorithm detailed in Subsection \ref{subse:existingAlgs}.
 
To ensure that all edges have the opportunity to be repositioned and to increase the likelihood of finding an improved solution, the general shift move is performed $k_{max}^3$ times during this procedure.

\begin{algorithm}[tbh]
\small
\caption{Local search procedure}
\label{alg:m-shift}
\SetKwInOut{Input}{Input}
\SetKwInOut{Output}{Output}
\Input{Current sequence $ S= (e_1, \ldots, e_m)$, number of attempts $k_{max}^3$, a positive integer $\theta$}
\Output{Updated sequence $S$}
\For{$i = 1, \ldots, k_{max}^3$}
{
    Select \( \theta \) random edges from \( S \) and assign them to \( T \)\;
    Remove all edges in $T$ from the sequence $S$\;
    \ForEach{$e \in T$}
    {
        Insert $e$ into S at the cheapest position\;
        Remove $e$ from $T$\;
    }
}
\end{algorithm}

\begin{figure}[tbh]
    \centering
    \begin{minipage}{0.45\textwidth}
        \centering
        \includegraphics[width=\textwidth]{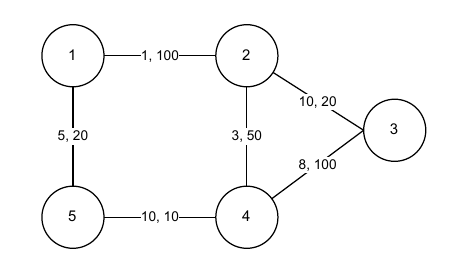}
        \caption{CPP-LC instance}
        \label{fig:m-shift-example}
    \end{minipage}%
    \hfill
    \begin{minipage}{0.45\textwidth}
        \centering
        \captionof{table}{Costs for inserting edge $\{2,4\}$.}
        \label{tab:m-shift-example-1}
        \begin{adjustbox}{width=1.0\textwidth}
        \begin{tabular}{l c}
            \toprule
            Sequence $S$ & Cost \\
            \midrule
            \textbf{\{2, 4\}}, \{1, 2\}, \{3, 2\}, \{5, 4\}, \{3, 4\} & 6455 \\
            \{1, 2\}, \textbf{\{2, 4\}}, \{3, 2\}, \{5, 4\}, \{3, 4\} & 5045 \\
            \{1, 2\}, \{3, 2\}, \textbf{\{2, 4\}}, \{5, 4\}, \{3, 4\} & 5725 \\
            \{1, 2\}, \{3, 2\}, \{5, 4\}, \textbf{\{2, 4\}}, \{3, 4\} & 6435 \\
            \{1, 2\}, \{3, 2\}, \{5, 4\}, \{3, 4\}, \textbf{\{2, 4\}} & 7315 \\
            \bottomrule
        \end{tabular}
        \end{adjustbox}

        \vspace{0.5cm} 

        \captionof{table}{Costs for inserting edge $(1,5)$.}
        \label{tab:m-shift-example-2}
        \begin{adjustbox}{width=1.0\textwidth}
        \begin{tabular}{l c}
            \toprule
            Sequence $S$ & Cost \\
            \midrule
  \textbf{\{1, 5\}}, \{1, 2\}, \{2, 4\}, \{3, 2\}, \{5, 4\}, \{3, 4\} & 7895 \\
  \{1, 2\}, \textbf{\{1, 5\}}, \{2, 4\}, \{3, 2\}, \{5, 4\}, \{3, 4\} & 7295 \\
  \{1, 2\}, \{2, 4\}, \textbf{\{1, 5\}}, \{3, 2\}, \{5, 4\}, \{3, 4\} & 7265 \\
  \{1, 2\}, \{2, 4\}, \{3, 2\}, \textbf{\{1, 5\}}, \{5, 4\}, \{3, 4\} & 5555 \\
  \{1, 2\}, \{2, 4\}, \{3, 2\}, \{5, 4\}, \textbf{\{1, 5\}}, \{3, 4\} & 6265 \\
  \{1, 2\}, \{2, 4\}, \{3, 2\}, \{5, 4\}, \{3, 4\}, \textbf{\{1, 5\}} & 6235 \\
            \bottomrule
        \end{tabular}
        \end{adjustbox}
    \end{minipage}
\end{figure}

To illustrate the general shift move, consider the graph $G$ depicted in Figure~\ref{fig:m-shift-example}, where the curb weight of the vehicle is $W=0$ and $\theta = 2$. The two numbers shown at an edge $e$ correspond to $d_e$ and $q_e$. Assume that the initial sequence is $S = (\{1, 2\}, \{1, 5\}, \{3, 2\}, \{2, 4\}, \{5, 4\}, \{3, 4\})$ and $T=\{\{2, 4\}, \{1, 5\}\}$. After removing the edges in $T$ from $S$, the current sequence becomes $S = (\{1, 2\}, \{3, 2\}, \{5, 4\}, \{3, 4\})$.

Table~\ref{tab:m-shift-example-1} presents the cost of all sequences obtained by inserting $\{2, 4\}$ into every position of $S$. The insertion position that yields the lowest cost, $5045$, corresponds to position 2. Thus, the edge $\{2, 4\}$ is inserted at position 2, and the updated sequence is $S = (\{1, 2\}, \{2, 4\}, \{3, 2\}, \{5, 4\}, \{3, 4\})$. 

The general shift move then proceeds to insert the remaining edge into $S$. Table~\ref{tab:m-shift-example-2} shows the same for $\{1, 5\}$. One can check that the final sequence obtained after this move is $S = (\{1, 2\}, \{2, 4\}, \{3, 2\}, \{1, 5\}, \{5, 4\}, \{3, 4\})$ and the total cost has been reduced from $7975$ to $5555$.

Evaluating each position to reinsert each edge in $T$ requires $O(m)$ time, which leads to a computational cost of $O(m^2)$ to determine the cheapest position for that edge. Consequently, each general shift move is computed in $O(\theta m^2)$ time, and the overall time complexity of the local search procedure is $O(k_{max}^3 \theta m^2)$.

\subsection{Reordering procedure}

After applying the local search, the resulting route may include edges that are deadheaded before being serviced, which violates the constraint~\ref{eq:constraint9} in the MILP model. To ensure that the route is ready for the intensification phase, a reordering phase is applied. This phase rearranges the sequence so that each edge is serviced the first time it is visited. After the sequence is adjusted, the deadheading distances between consecutive edges are recalculated to account for the new order of servicing.

These adjustments offer an additional advantage, as they enable the vehicle to reduce its load earlier in the route, thereby reducing subsequent deadheading costs and, consequently, the overall cost. This will become clear in the following example. Figure~\ref{fig:reordering} presents an illustrative example of the reordering procedure applied to the best sequence obtained from Table~\ref{tab:m-shift-example-2}. The first row shows the current sequence obtained after the general shift move. The second row shows the actual route corresponding to this sequence after applying dynamic programming. For clarity, edges that are serviced are shown as solid lines, while deadheaded edges are represented by dashed lines.

\begin{figure}[htb]
    \centering
    \begin{adjustbox}{width=\textwidth}
    \begin{tikzpicture}[>=stealth, node distance=1.5cm]
        \tikzstyle{dot} = [circle, fill=black, inner sep=1.5pt]

        \node at (-2, 6) {\Large \textbf{Before reordering:}};
        
        \node at (-2, 4.5) {\Large \textbf{Sequence:}};
        \node at (4.2, 4.5) {\Large (\{1,2\}, \{2,4\}, \{3,2\}, \{1,5\}, \{5,4\}, \{3,4\})};

        \node at (-2.4, 3) {\Large \textbf{Route:}};
        \node[dot] (A1) at (0,3) [label=below:\Large 1] {};
        \node[dot] (B1) at (1.5,3) [label=below:\Large 2] {};
        \node[dot] (C1) at (3,3) [label=below:\Large 4] {};
        \node[dot] (D1) at (4.5,3) [label=below:\Large 3] {};
        \node[dot] (E1) at (6,3) [label=below:\Large 2] {};
        \node[dot] (F1) at (7.5,3) [label=below:\Large 1] {};
        \node[dot] (G1) at (9,3)  [label=below:\Large 5] {};
        \node[dot] (H1) at (10.5,3) [label=below:\Large 4] {};
        \node[dot] (I1) at (12,3) [label=below:\Large 3] {};
        \node[dot] (J1) at (13.5,3) [label=below:\Large 2] {};
        \node[dot] (K1) at (15,3) [label=below:\Large 1] {};

        \draw[thick,-] (A1) -- (B1);
        \draw[thick,-] (B1) -- (C1);
        \draw[dashed,-] (C1) -- (D1);
        \draw[thick,-] (D1) -- (E1);
        \draw[dashed,-] (E1) -- (F1);
        \draw[thick,-] (F1) -- (G1);
        \draw[thick,-] (G1) -- (H1);
        \draw[thick,-] (H1) -- (I1);
        \draw[dashed,-] (I1) -- (J1);
        \draw[dashed,-] (J1) -- (K1);

        \node at (-2, 1) {\Large \textbf{After reordering:}};
        
        \node at (-2, -0.5) {\Large \textbf{Sequence:}};
        \node at (4.2,-0.5) {\Large(\{1,2\}, \{2,4\}, \{3,4\}, \{3,2\}, \{1,5\}, \{5,4\})};

        \node at (-2.4, -2) {\Large \textbf{Route:}};
        \node[dot] (A3) at (0,-2) [label=below: \Large 1] {}; 
        \node[dot] (B3) at (1.5,-2) [label=below:\Large 2] {};
        \node[dot] (C3) at (3,-2) [label=below:\Large 4] {};
        \node[dot] (D3) at (4.5,-2) [label=below:\Large 3] {};
        \node[dot] (E3) at (6,-2) [label=below:\Large 2] {};
        \node[dot] (F3) at (7.5,-2) [label=below:\Large 1] {};
        \node[dot] (G3) at (9,-2) [label=below:\Large 5] {};
        \node[dot] (H3) at (10.5,-2) [label=below:\Large 4] {};
        \node[dot] (I3) at (12,-2) [label=below:\Large 2] {};
        \node[dot] (J3) at (13.5,-2) [label=below:\Large 1] {};

        \draw[thick,-] (A3) -- (B3);
        \draw[thick,-] (B3) -- (C3);
        \draw[thick,-] (C3) -- (D3);
        \draw[thick,-] (D3) -- (E3);
        \draw[dashed,-] (E3) -- (F3);
        \draw[thick,-] (F3) -- (G3);
        \draw[thick,-] (G3) -- (H3);
        \draw[dashed,-] (H3) -- (I3);
        \draw[dashed,-] (I3) -- (J3);
    \end{tikzpicture}
    \end{adjustbox}
    \caption{An example of the Reordering procedure.}
    \label{fig:reordering}
\end{figure}

In this route, one can see that edge $\{3, 4\}$ is serviced in position 6 in the current sequence $S$, but it is deadheaded between edges $\{2, 4\}$ and $\{3, 2\}$ in positions 2 and 3, respectively. The reordering phase moves the edge $\{3, 4\}$ to position 3 in $S$ and shortens the deadheading distance from the edges $\{5, 4\}$ to the depot. The new sequence and its corresponding route are shown in the third and fourth rows, respectively. We can check that the total cost is reduced from $5,555$ to $2,155$.

To calculate the route for the current sequence, we use dynamic programming in $O(m)$ time to determine the traversal direction of all edges. It should be noted that the shortest paths between any two nodes in the graph have been precomputed, so checking the deadheading status of edges in the shortest path between two consecutive edges can be done in $O(m)$ time. Therefore, the overall time complexity of this procedure is $O(m^2)$. Details of this reordering phase are provided in Algorithm~\ref{alg:preprocess}.

\begin{algorithm}[tbh]
\small
\caption{Reordering procedure}
\label{alg:preprocess}
\SetKwInOut{Input}{Input}
\SetKwInOut{Output}{Output}
\Input{Sequence \( S= (e_1, \ldots, e_m) \)}
\Output{Updated sequence \( S \)}
Run dynamic programming to find the route corresponding to $S$\;
Let $i=1$\;
\While{$i<m-1$}
{
    Let \( P \) be the shortest path from the end node of edge \( e_i \) to the start node of edge \( e_{i+1} \)\; 
    Let $j=i+1$\;
    \ForEach{$e \in P$}
    {
        \If{the position of $e$ in $S$ is greater than $i+1$}
        {
            Move the edge $e$ to position where it is deadheaded\;
            Increase $j$ by 1\;
        }
    }
    Let $i=j$\;
}
\end{algorithm}

\subsection{Intensification procedure}
The intensification procedure begins once the reordering procedure is completed. It deeply explores the specific neighborhood of the current best solution for potential improvements. This neighborhood is defined by selecting a subsequence of consecutive edges from the current sequence and allowing their permutation within the subsequence, while keeping the remaining edges fixed.

To incorporate this neighborhood into the MILP formulation, we first convert the current route to obtain the MILP solution $(\bar{x}, \bar{y}, \bar{f})$. This conversion can be performed greedily in linear time, with details provided in Algorithm~\ref{alg:conversion}.

Next, we select a subsequence of $\delta$ consecutive edges, denoted as $(e_{\epsilon}, e_{\epsilon+1}, \dots, e_{\epsilon+\delta-1})$, where $\epsilon$ is the starting position ($1 \leq \epsilon \leq m-\delta+1$). Let $i_k$ and $j_k$ denote the start and end nodes, respectively, of edge $e_k$, following its traversal direction in the current route. We then fix certain variables in the MILP as follows:
\begin{align}
     y_{i_k, j_k}^k &= \bar{y}_{i_k, j_k}^k & \forall k \in \{1, \dots, \epsilon-1\} \cup \{\epsilon+\delta, \dots, m\} \label{eq:y-fix}\\
     f_k &= \bar{f}_k & \forall k \in \{1, \dots, \epsilon-1\} \cup \{\epsilon+\delta, \dots, m\} \label{eq:f-fix}\\
     x_{i_k, j_k}^k &= \bar{x}_{i_k, j_k}^k & \forall k \in \{1, \dots, \epsilon-2\} \cup \{\epsilon+\delta, \dots, m\}\label{eq:x-front-fix}
\end{align}
The equalities \ref{eq:y-fix} fix the required edges in unchosen segments (or positions), while the equalities \ref{eq:f-fix} enforce that the vehicle load during these unchosen segments remains unchanged. The variables $x$, which determine the deadheaded edges for each segment, are also fixed for the unchosen segments, except for period $\epsilon-1$. This helps improve deadheading costs between the required edges at positions $\epsilon-1$ and $\epsilon$.

Note that the number of free variables in the reduced MILP is $O(m\delta)$, while in the original formulation, it is $O(m^2)$. This reduced MILP is significantly smaller than the original and much easier to solve when $\delta$ is small enough. The choice of $\delta$ will be based on our preliminary experiments and is provided in Section \ref{sec:experiments}. To ensure that every edge has an opportunity to be moved, $\epsilon$ can be selected randomly each time we call the intensification procedure.

\begin{algorithm}[tbh]
\small
\caption{Transforming current route into MILP solution}
\label{alg:conversion}
\SetKwInOut{Input}{Input}
\SetKwInOut{Output}{Output}
\Input{Current sequence \( S= (e_1, \ldots, e_m) \) and its corresponding route}
\Output{MILP solution \( (\bar{x}, \bar{y}, \bar{f}) \)}

Initialise vectors \( \bar{x}, \bar{y}, \bar{f} \) as null vectors\; 
Set $\bar{f}_{1} = Q$ and $\bar{f}_{m+1} =0$\; 
\For{$k = 1, \ldots, m-1$}
{
    \tcc{edge \( e_k =\{i_k, j_k\}\) is served from $i_k$ to $j_k$ in period $k$}
    Set \( \bar{y}^{k}_{i_k, j_k} = 1 \)\;
    Set \( \bar{f}_{k+1} \) to $\bar{f}_{k} - q_{e_k}$ \;
    Let \( P \) be the deadheading path from \( j_k \) to \( i_{k+1} \) \; 
    
    \ForEach{$e=\{i,j\} \in P$}
    {
         \tcc{edge \( e \) is deadheaded from $i$ to $j$ in period $k$}
        Set \( \bar{x}^k_{ij} = 1 \)\;
    }
}
    \tcc{edge \( e_m \) is served from $i_m$ to $j_m$ in the last period}
    Set \( \bar{y}^{m}_{i_m, j_m} = 1 \)\;
    Let \( P \) be the deadheading path from the end node of edge \( e_m \) to depot\;
     \ForEach{$e=\{i,j\} \in P$}
    {
         \tcc{edge \( e \) is deadheaded from $i$ to $j$ in the last period}
        Set \( \bar{x}^m_{ij} = 1 \)\;
    }
\end{algorithm}

\section{Ant Colony Optimization}
\label{sec:aco}

Cober\'an et al. \cite{corberan2018} introduced two well-known metaheuristics, ILS and VNS, for the CPP-LC. To extend this study, we introduce the ACO metaheuristic, which is widely applied in the vehicle routing literature, to CPP-LC. Our objective is to evaluate whether ACO can effectively address the load-dependent cost feature and to compare its performance with that of ILS and VNS. Since both ILS and VNS were implemented in a simplified form, we also design a basic version of ACO to ensure a fair comparison.

The ACO algorithm is inspired by the foraging behavior of real ants, a phenomenon first studied by biologist Grass\'e in 1959 \cite{grasse1959reconstruction}. In nature, ants initially search for food at random and deposit small amounts of a chemical substance known as pheromone along their paths. After locating a food source, the ants reinforce the trail by depositing an additional pheromone. The quality of the food source influences the amount of pheromone deposited, with better sources receiving more.  Other ants are more likely to follow paths with higher levels of pheromones, and if they also find food, they further strengthen the trail. Less favorable paths gradually lose pheromone through evaporation and eventually disappear. Over time, this feedback mechanism leads the colony to converge on the most promising routes. The ACO variants differ in how the ants select paths and how pheromone levels are updated. A comprehensive overview of ACO algorithms and their applications can be found in, e.g., \cite{dorigo2004ant, dorigo2018ant}.

We define the state space $\mathcal{S}$, where each state is represented by a tuple $(i, j, d)$. Here, $(i, j) \in E$ denotes an edge that requires servicing, and $d \in \{1, 2\}$ indicates the direction of service: $d = 1$ means the edge is serviced from $i$ to $j$, and $d = 2$ represents the reverse direction from $j$ to $i$. In addition, an auxiliary state $s^*$ is introduced to indicate that the ant is at the depot and has not yet visited any other states. Therefore, the total number of states is $2m + 1$.

Each ant constructs a solution by traversing a sequence of $m$ states. At each step of the algorithm, an ant transitions from state $s$ to state $t$ with a probability $p(t|s)$, which is computed using two sources of information: the attractiveness of the move $\eta_{st}$, representing the a priori desirability of the move, and the level of pheromones $\tau_{st}$, representing the a posteriori desirability. The transition probability is given by:

\begin{equation}\label{eq:transitionProb}
    p(t|s) = \frac{\tau_{st} \cdot \eta_{st}}{\sum_{t' \in U(s)} \tau_{st'} \cdot \eta_{st'}},
\end{equation}

where $U(s)$ denotes the set of unvisited states.

The attractiveness $\eta_{st}$ of a move from state $s = (i,j, d)$ to state $t = (i', j', d')$ is computed as follows:
\begin{equation}
    \eta_{st} = \sqrt{(W + q_{e'}) D_{s, t} + \left(W + \frac{q_{e'}}{2}\right) d_{e'}},
    \label{eq:init-eta-1}
\end{equation}
where $D_{s, t}$ is the shortest distance from the end node of edge $e = \{i,j\}$ in direction $d$ to the start node of edge $e'=\{i',j'\}$ in direction $d'$. In the case where $x = s^*$, we use the following equation instead:

\begin{equation}
    \eta_{s^*t} = \sqrt{(W + Q) D_{1, t} + \left(W + Q - \frac{q_{e'}}{2}\right) d_{e'}},
    \label{eq:init-eta-2}
\end{equation}
where $D_{1, t}$ is the shortest distance from the depot to the start node of edge $e'$ in direction $d'$.

After all ants have completed their paths, the pheromone levels are updated. Based on the quality of the solutions, the pheromone levels are adjusted accordingly, either reinforced or evaporated, as follows:
\begin{equation}
\tau_{st} \leftarrow (1 - \rho) \tau_{st} + \sum_{a = 1}^{p_{\text{max}}} \Delta \tau^a_{st},
\label{eq:pheromone-update}
\end{equation}

where $\rho$ is the coefficient of evaporation of pheromones, $p_{\text{max}}$ is the number of ants and $\Delta \tau^a_{st}$ represents the amount of pheromone deposited by the $a$-th ant, which is defined as:

\begin{equation}
\Delta \tau^a_{st} = 
\begin{cases}
\frac{C}{\sqrt{z(S_a)}} & \text{if the ant transitions from state $x$ to $y$}, \\
0 & \text{otherwise}.
\end{cases}
\end{equation}

Here, $S_a$ represents the sequence of states visited by the $a$-th ant, $z(S_a)$ is the total cost of the sequence $S_a$, and $C$ is a constant.

Algorithm~\ref{algo:ACO} presents the pseudocode for the ACO algorithm, while Algorithm~\ref{algo:sampling} presents the sampling procedure used in ACO. To accelerate the sampling process, we employed parallelism by executing lines 7 to 10 concurrently across multiple threads.
\begin{algorithm}[tbh]
\small
 \caption{ACO for CPP-LC}
    \label{algo:ACO}
\SetKwInOut{Input}{Input}
\SetKwInOut{Output}{Output}
\Input{Graph $G = (V, E)$, evaporation constant $\rho = 0.8$,  $\epsilon = 0.001$,
\\ $C = 1$, maximum number of iterations $k_{\text{max}}^4$, number of artificial ants $p_{\text{max}}$}
\Output{The best sequence of states $S^*$}
    $\tau_{st} \leftarrow \epsilon \ \forall s, t \in \mathcal{S}$\;
    Initialize $\eta_{st}$ by Eq.~(\ref{eq:init-eta-1}) and Eq.~(\ref{eq:init-eta-2})\;
    \For{$k = 1 \rightarrow k_{\text{max}}^4$}{
        \For{$a = 1 \rightarrow p_{\text{max}}$}{
            Sample $S_a$ by Algorithm~\ref{algo:sampling}\;
            \If{$z(S_a) < z(S^*)$}{
                $S^* \leftarrow S_a$\;
            }
        }
        Update the pheromone by Eq.~(\ref{eq:pheromone-update})\;
    }
\end{algorithm}

\begin{algorithm}[tbh]
\SetKwInOut{Input}{Input}
\SetKwInOut{Output}{Output}
\Input{$\tau_{st}$, $\eta_{st}$, and the starting state $s^*$}
\Output{The sample sequence of states $S$}
    Initialise the set of unvisited edges $\bar{E} \leftarrow E$ \;
    $s \leftarrow s^*$ \;
    $S \leftarrow \emptyset$ \;
    \For{$i = 1 \rightarrow |E|$}{
        Construct the set of unvisited states $U(s) = \{(i, j, d) \mid (i, j) \in \bar{E}, d \in \{1,2\}\}$ \;
        Compute the transition probabilities using Eq. (\ref{eq:transitionProb})\;
        Sample $s'$ from the transition probabilities. Assume $s' = (i, j, d)$ \;
        $s \leftarrow s'$ \;
        $S \leftarrow S \cup (i, j, d)$ \;
        $\bar{E} \leftarrow \bar{E} \setminus (i, j)$\;
    }
    \caption{Sampling procedure of ACO for CPP-LC}
    \label{algo:sampling}
\end{algorithm}

\section{Computational experiments}
\label{sec:experiments}

\subsection{Experimental setup}

All experiments were performed on a system equipped with an AMD Ryzen 5 5600X 6-core processor (3.70 GHz) and 32 GB of RAM, running Windows 11 with the MSVC compiler (version 19.29). We used the MILP solver of CPLEX 12.10 with its default settings to solve the reduced MILPs. A time limit of 100 seconds per reduced MILPs was imposed.

The performance of the proposed algorithms is evaluated on three distinct sets of CPP-LC instances generated following the procedure described in \cite{corberan2018}. Edge demands $q_e$ are generated under two different scenarios: (i) a \emph{proportional demand} setting, where $q_e = d_e$ for all $e \in E$, and (ii) a \emph{non-proportional demand} setting, where $q_e$ values are randomly generated. These two scenarios will allow us to assess the robustness of the algorithm in both structured and randomly varying demand patterns. For each instance, three values of the curb weight parameter $W$ are considered: $W = 0$, $W = \frac{Q}{2}$, and $W = 5Q$, where the total edge demand is defined as $Q = \sum_{e \in E} q_e$. The three sets of CPP-LC instances are as follows.

\begin{itemize}
    \item \textbf{First dataset:} This dataset includes 18 small Eulerian instances with $|V| \in \{7, 8\}$ and $|E| \in \{8, 9, 10\}$. For these instances, their optimal route can be determined by brute-force search. 
     \item \textbf{Second dataset:} This dataset comprises 60 small-to-medium-sized instances, including 18 Eulerian graphs, 18 graphs from Christofides et al. \cite{christofides1981algorithm}, and 24 graphs from Hertz et al. \cite{hertz1999improvement}. The number of edges $|E|$ ranges from 11 to 33. This dataset is used to evaluate the performance of the algorithms across diverse graph structures and demand configurations.
    \item \textbf{Third dataset:} This dataset contains 48 larger Eulerian instances. The instances are constructed with $|V| \in \{10, 20, 30\}$ and $|E| \in \{16, 20, 27, 53, 75, 110, 162, 232\}$. This dataset is used to test the scalability and robustness of algorithms in larger instances.
\end{itemize}

For all experiments, we set the number of outer iterations at $k_{\max}^1 = 5$, and the maximum number of intensification attempts per outer iteration in the MaLD scheme at $k_{\max}^2 = 2$. We also set the parameter $\theta$ to $m/2$ in the local search phase. In ACO, we set the number of iterations $k_{\max}^4$ to $1000$ and the number of artificial ants $p_{\max}$ to $5$. For comparison, we also implement ILS and VNS metaheuristics in \cite{corberan2018} with parameter $k_{\max} = 75$. 

To handle varying instance sizes, the following parameters are adjusted: 
\begin{itemize}
    \item \textbf{First and second datasets}: The parameter $\delta$ in the intensification phase is set to $\frac{m}{2}$, and the number of attempts $k_{\text{max}}^3$ in the local search procedure is set to $100$.
    \item \textbf{Third dataset}: The parameter $\delta$ in the intensification phase is reduced to $\frac{m}{4}$, and the number of attempts $k_{\text{max}}^3$ in the local search procedure is reduced to $35$.
\end{itemize}

These parameters were determined through extensive preliminary experiments performed on modified datasets, in which data values were randomly varied by up to $10\%$.
\subsection{Experimental results}

Tables~\ref{tab:res_1}–\ref{tab:res_3} present the summary results of five algorithms: GCH, ILS, VNS, ACO, and MaLD, on different datasets. The instances in each dataset are grouped according to different criteria: by instance source in Table~\ref{tab:res_2}, and by edge set size in Tables~\ref{tab:res_1} and~\ref{tab:res_3}. Each table reports the average percentage gap with respect to the best-known solution and the average running time (in seconds) for each algorithm. 

\begin{table}[tbh]
  \centering
  \setlength\tabcolsep{2pt}
    \setlength\extrarowheight{2pt} 
  \caption{Performance of GCH, ILS, VNS, ACO, and MaLD on the first dataset of small Eulerian CPP-LC instances with $|E| \in \{8,9,10\}$. Reported values are the average percentage gap from the optimal solution and the average running time (in seconds).}
  \begin{tabular*}{\textwidth}{@{\extracolsep{\fill}}{l}@{\extracolsep{\fill}}*{10}{r}}
    \toprule
    && \multicolumn{2}{c}{$|E|=8$}    && \multicolumn{2}{c}{$|E|=9$}   && \multicolumn{2}{c}{$|E|=10$}   \\[0.1cm]
    \cline{3-4} \cline{6-7} \cline{9-10}
   Algorithms &&\textbf{Gap (\%)} & \textbf{Time (s)}&&\textbf{Gap (\%)} & \textbf{Time (s)}&&\textbf{Gap (\%)} & \textbf{Time (s)}\\
    \midrule
GCH	&&	1.12	&	0.000	&&	1.85	&	0.000	&&	1.54	&	0.001	\\
ILS	&&	0.00	&	0.076	&&	0.00	&	0.110	&&	0.00	&	0.140	\\
VNS	&&	0.24	&	0.048	&&	0.00	&	0.069	&&	0.00	&	0.086	\\
ACO	&&	1.98	&	0.479	&&	4.88	&	0.483	&&	1.59	&	0.482	\\
MaLD	&&	\textbf{0.00}	&	0.274	&&	\textbf{0.00}	&	0.382	&&	\textbf{0.00}	&	0.536	\\
    \bottomrule
  \end{tabular*}
  \label{tab:res_1}
\end{table}

\begin{table}[tbh]
  \centering
  \setlength\tabcolsep{2pt}
    \setlength\extrarowheight{2pt} 
  \caption{Performance of GCH, ILS, VNS, ACO, and MaLD on the second dataset of small-to-medium CPP-LC benchmark instances, grouped by instance source. Reported values are the average percentage gap from the best-known solution and the average running time (in seconds).}
  \begin{tabular*}{\textwidth}{@{\extracolsep{\fill}}{l}@{\extracolsep{\fill}}*{10}{r}}
    \toprule
    && \multicolumn{2}{c}{\textbf{Eulerian}}    && \multicolumn{2}{c}{\textbf{Christofides et al.}}   && \multicolumn{2}{c}{\textbf{Hertz et al.}}   \\[0.1cm]
    \cline{3-4} \cline{6-7} \cline{9-10}
   Algorithms &&\textbf{Gap (\%)} & \textbf{Time (s)}&&\textbf{Gap (\%)} & \textbf{Time (s)}&&\textbf{Gap (\%)} & \textbf{Time (s)}\\
    \midrule
GCH && 17.61 & 0.007 && 7.69 & 0.002 && 11.44 & 0.002 \\
ILS && 8.18 & 1.384 && 4.05 & 0.406 && 4.38 & 0.450 \\
VNS && 8.46 & 0.668 && 4.48 & 0.189 && 5.53 & 0.215 \\
ACO && 10.51 & 0.693 && 9.90 & 0.540 && 9.08 & 0.544 \\
MaLD && \textbf{0.00} & 13.372 && \textbf{0.00} & 110.230 && \textbf{0.00} & 164.919 \\
    \bottomrule
  \end{tabular*}
  \label{tab:res_2}
\end{table}

\begin{table}[tbh]
  \centering
  \setlength\tabcolsep{2pt}
    \setlength\extrarowheight{2pt} 
  \caption{Performance of GCH, ILS, VNS, ACO, and MaLD on the third dataset of large CPP-LC instances, grouped by edge-set size. Reported values are the average percentage gap from the best-known solution and the average running time (in seconds).}
  \begin{tabular*}{\textwidth}{@{\extracolsep{\fill}}{l}@{\extracolsep{\fill}}*{10}{r}}
    \toprule
    && \multicolumn{2}{c}{$|E|\in \{16, 20, 27\}$}    && \multicolumn{2}{c}{$|E|\in \{53,75\}$}   && \multicolumn{2}{c}{$|E|\in \{110, 162, 232\}$}   \\[0.1cm]
    \cline{3-4} \cline{6-7} \cline{9-10}
   Algorithms &&\textbf{Gap (\%)} & \textbf{Time (s)}&&\textbf{Gap (\%)} & \textbf{Time (s)}&&\textbf{Gap (\%)} & \textbf{Time (s)}\\
    \midrule
   GCH	&&	8.31	&	0.006	&&	11.08	&	0.120	&&	2.25	&	2.260	\\
ILS	&&	6.05	&	1.131	&&	6.21	&	27.181	&&	1.35	&	531.647	\\
VNS	&&	6.23	&	0.528	&&	6.21	&	13.195	&&	1.35	&	259.085	\\
ACO	&&	6.98	&	0.611	&&	6.21	&	2.089	&&	1.35	&	14.203	\\
MaLD	&&	\textbf{0.00}	&	8.453	&&	\textbf{0.00}	&	553.911	&&	\textbf{0.00}	&	1262.218	\\
    \bottomrule
  \end{tabular*}
  \label{tab:res_3}
\end{table}

Table~\ref{tab:res_1} reports the results on the first dataset, which includes the smallest instances with $|E| \in \{8, 9, 10\}$, for which the optimal solutions are known. Among the algorithms, GCH is the fastest and, unexpectedly, outperforms ACO in both solution quality and computation time. The possible reasons for this are as follows: (1) in small-scale instances, ACO incurs significant computational overhead due to the involvement of multiple agents, pheromone updates, and probabilistic path construction; (2) ACO's random exploration often wastes time on less promising areas, delaying the identification of high-quality solutions; and (3) the pheromone learning mechanism may become ineffective in small solution spaces.

The remaining algorithms successfully improve the initial solutions produced by GCH within fractions of a second. In particular, MaLD consistently found the optimal solution in every instance, with an average computation time of $0.4$ seconds. ILS and VNS were even faster, requiring less than $0.15$ seconds on average, and both reached optimality in nearly all cases. Given the small size of these instances, differences in computation time between the algorithms are negligible and likely due to system-level noise.

Table~\ref{tab:res_2} reports the results for the second dataset, which consists of small to medium instances with $|E|$ ranging from $11$ to $33$. GCH remains the fastest, but less accurate. For Eulerian and Hertz et al. instances, ACO outperforms GCH in terms of solution quality, although it incurs significantly higher computational costs. In contrast, for the Christofides et al. instances, ACO performs worse than GCH in both solution quality and running time.

MaLD, ILS, and VNS further improve upon the solutions provided by GCH. However, the solutions provided by ILS and VNS are 4-9\% worse than those obtained by MaLD. However, MaLD requires substantially more computational time, particularly for Hertz et al. instances.

Table~\ref{tab:res_3} presents the results for the largest dataset, with instances grouped by edge set size as follows: $|E| \in \{16, 20, 27\}$, $|E| \in \{53, 75\}$, and $|E| \in \{110, 162, 232\}$. These groupings are based on the similarity of the results across algorithms. In this dataset, ACO consistently outperforms GCH in terms of solution quality, while being the second-fastest algorithm after GCH. Furthermore, ACO achieves solution quality comparable to that of ILS and VNS, but with significantly lower computational time. This demonstrates the effectiveness of ACO's exploratory capabilities and its learning mechanism in larger instances. While ILS and VNS exhibit substantial increases in running time due to the increasing number of neighbor evaluations, ACO's compatibility with parallel execution enables it to scale more efficiently in larger instances.

It should be noted that for instances where $|E| \ge 53$, the metaheuristic approaches, ILS, VNS, and ACO, produce similar local optima, while the matheuristic method, MaLD, is able to escape these local minima. Specifically, MaLD improves the solutions obtained by ILS and VNS by approximately 6 to 7\% for instances where $16 \leq |E| \leq 75$ and by an average of 1.35\% for instances with $|E| \in \{110, 162, 232\}$. The modest improvement in the latter group is because we adjusted the parameter $\delta$ from $m/2$ to $m/4$ to strike a better balance between the quality of the solution and the computational time required to solve the reduced MILPs.

Although MaLD remains the most computationally intensive method overall, we tested increasing the maximum number of attempts, $k_{\text{max}}$, to 150 and 300 for ILS and VNS. However, this adjustment did not improve solution quality but resulted in longer run times than MaLD. For brevity, these results have been omitted.

\subsection{Effect of MaLD components}
To assess the individual contributions of two key procedures, namely intensification and local search, to the overall performance of MaLD, we performed an experiment in which each procedure was removed and the results are summarized in Table~\ref{tab:ablation_summary}. The first column lists the instance groups, which are defined either by the size of the edge set or by the source of the benchmark instances, as in the previous experiment. The second and third columns report the average objective value and the average running time in seconds, respectively, achieved by the full MaLD algorithm. The last two columns show the average deviations in objective value and computational time when either the intensification or the local search procedure is disabled. Positive deviations indicate a deterioration in performance, such as an increase in objective cost or a longer running time, while negative deviations indicate an improvement.

Across most instance groups, removing either procedure results in a degradation of solution quality, with a greater impact observed in larger instances. For example, in the $|E| \in \{110, 162, 232\}$ group, excluding the intensification procedure increases the objective value by an average of 211,696, compared to 119,391 when omitting the local search procedure. A similar pattern is observed in the $|E|\in \{53, 75\}$ group and in the instances of Eulerian and Hertz et al., indicating that intensification plays a more critical role in improving solution quality in these instances. However, this pattern is not observed in cases with $|E| \in \{16, 20, 27\}$, where both components contribute significantly to the quality of the solution, with the local search having a slightly greater impact.

In small instances with $|E| \in \{8, 9, 10\}$, removing either component has negligible impact on solution quality. This suggests that since GCH already produces solutions within 2\% of optimality and the search space is small, either procedure alone is often sufficient to obtain high-quality solutions.

An exception occurs in the case of Christofides et al., where the use of only intensification slightly improves the quality of the solution ($-3$) but increases runtime. In contrast, excluding intensification leads to worse outcomes ($+350$).

In terms of running time, the reductions resulting from component removal are generally substantial. On average, excluding the intensification procedure reduces the runtime by approximately 74\%, while removing the local search component results in a reduction of approximately 20\%.


\begin{table}[htbp]
\centering
\setlength\tabcolsep{2pt}
    \setlength\extrarowheight{2pt} 
\caption{Impact of MaLD component removal on objective value and runtime}
\label{tab:ablation_summary}
\small
\begin{tabular*}{\textwidth}{@{\extracolsep{\fill}}{l}@{\extracolsep{\fill}}*{10}{r}}
\toprule
&& \multicolumn{2}{c}{Full MaLD} &&\multicolumn{2}{c}{Deviation (- Intensification / - LocalSearch)}\\
\cline{3-4} \cline{6-7}
\textbf{Instance groups} && \textbf{Obj. Cost} & \textbf{Time (s)} && \textbf{Obj. Cost} & \textbf{Time (s)}\\
\midrule
$|E|\in\{110,162,232\}$ && 11,857,017 & 1262.22 && 211,696 / 119,391 & - 818.66 / - 397.09 \\
$|E|\in \{53,75\}$ && 1,641,706 & 553.91 && 68,878 / 11,492 & - 467.60 / - 1.87 \\
$|E|\in \{16,20,27\}$ && 184,716.20 & 8.45 && 847.50 / 981.50 & - 5.60 /  4.99 \\
$|E|=8$ && 32,046 & 0.27 && 0 / 0 & - 0.08 / - 0.03 \\
$|E|=9$ && 23,737 & 0.38 && 210 / 210 & - 0.13 / - 0.07 \\
$|E|=10$ && 33,210 & 0.54 && 0 / 0 & - 0.18 / - 0.13 \\
Eulerian && 188,710 & 13.37 && 1,062 / 163 & - 5.44 / - 6.19 \\
Christofides et al. && 82,841 & 110.23 && 350 / -3 & - 107.72 /  6.67 \\
Hertz et al.&& 41,717 & 164.92 && 365 / 133 & - 161.89 / - 12.51 \\
\bottomrule
\end{tabular*}
\label{tab:matheuristic_fullsummary}
\end{table}

\section{Conclusion}
\label{sec:conclusion}

In this paper, we introduced a matheuristic for the CPP-LC. The proposed MaLD consistently achieved the best solution quality compared to existing algorithms, reaching optimality in all instances where the optimal solutions are known. The effectiveness of MaLD stems from its combination of local search, which efficiently explores promising regions of the solution space, and intensification, which uses reduced MILP models to deepen the search.

Although ILS and VNS have previously been introduced to CPP-LC, we introduced ACO to provide a comprehensive comparison. Experimental results show that ACO performs competitively in large instances because of its strong exploratory capabilities and compatibility with parallel execution. However, for smaller instances, its computational overhead and stochastic nature reduce its efficiency.

Future research may extend this study in several directions. First, a decomposition heuristic can be incorporated within the intensification phase to partition the solution space into independent subspaces. Each subspace can then be explored efficiently by solving reduced MILP models in parallel, thus improving computational time. Second, reinforcement learning techniques could be employed to automatically tune parameters within each component of the algorithm. We hope to continue in these directions in future work. Potential applications include smart logistics, waste collection, urban transportation planning, and energy-aware routing systems.

The implementation of our algorithms, instances, and detailed results are available at \url{https://github.com/HySonLab/MatCPP}.




\bibliographystyle{elsarticle-harv}
\bibliography{main}







\end{document}